\documentclass[12pt]{article}
\usepackage[fleqn]{amsmath}
\usepackage[cp1251]{inputenc}
\usepackage[russian,english]{babel}
\usepackage{amsmath,amsfonts,amssymb,graphicx,makeidx,amsthm,cite}
\usepackage{xcolor}
\usepackage{mathrsfs}
\newtheorem{teo}{Theorem}

\begin{document}
	
\begin{center}
	{\LARGE\bf Random walk and Weber-Schafheitlin integral: generalizations and discussion}
		\end{center}
		
\begin{center}
	{\LARGE\bf		Marichev O.I.$^{1}$, Shishkina E.L.$^{2,3,4,5}$ }
	\end{center}

E-mail: {\bf oleg@wolfram.com, ilina\_dico@mail.ru}

\vskip 0.5 cm

$^{1}$Wolfram Research, 100 Trade Center Dr, Champaign, IL 61820, USA

$^{2}$Department of Mathematical and Applied Analysis, Voronezh State University, Voronezh,  Universitetskaya pl. 1, 394018 Russia

$^{3}$Department of Applied Mathematics and Computer Modeling, Belgorod State National Research University (BelGU), Pobedy St., 85, Belgorod, 308015, Russia

$^{4}$Institute of Mathematics, Physics and Information Technology, Kadyrov Chechen State University, A. Sheripova st., 32, Grozny, 364024, Russia

$^{5}$International Laboratory of Stochastic Analysis and its Applications, National Research University Higher School of Economics, Pokrovsky Bulvar, 11, Moscow, 109028, Russia
 
\vskip 0.5 cm
	
{\bf Abstract.} This paper presents a study of the Weber–Schafheitlin integral and its connections to random walk theory. The Weber–Schafheitlin integral, defined as the improper integral of a product of Bessel functions, arises naturally in the evaluation of probability distributions for multidimensional random walks. This work bridges pure special function theory with applied probability.

	
	\section{Introduction}
	
	A random walk is a type of stochastic process (exemplified by Brownian motion) that models a trajectory composed of a sequence of random displacements within a given state space \cite{Verburgt}. 
	John Venn presented graph of a random walk and anticipated the notion of a fractal in the third edition of The Logic of Chance \cite{Venn}.
Karl Pearson \cite{PEARSON1,PEARSON2} described a walker who starts at a point, takes $n$ steps of fixed length in random directions, and asked for the probability of ending at a given distance from the origin. The path's aimless character, which Pearson himself noted, gave rise to the enduring nickname "the drunkard's walk."

The Weber-Schafheitlin integral is a specific type of improper integral that features a product of  Bessel functions in its integrand. 
Namely, for general values of the parameters and variables $\mu{\,\in\,}\mathbb{C}$, $\overline{\nu}{\,=\,}(\nu_1,...,\nu_m)$,
$\nu_j{\,\in\,}\mathbb{C}$, $\overline{x}=(x_1,...,x_m)$, $x_j{\,\in\,}\mathbb{R}$, $j{\,=\,}1,2,...,m$ involved, we will call the improper integral
\begin{equation}\label{WSh}
	\operatorname{WS}[\mu;\overline{\nu};\overline{x}] =\int\limits_{0}^{\infty} t^{\mu-1}\prod\limits_{j=1}^m J_{\nu_j}(x_jt)\,dt,
\end{equation}
   the {\bf Weber--Schafheitlin integral}, which we denote as $WS$.
Integral \eqref{WSh} was evaluated in \cite{Srivastava} in terms of the Lauricella hypergeometric function $F_C^{(m-1)}$ of $(m-1)$ variables.

 For the cases where $m{\,=\,}1$ the integral \eqref{WSh} is
known as the Mellin transform of the Bessel function of the first kind (see \cite{Hardy}, \cite{Luke}, p. 42).
When $m{\,=\,}2$, the expression \eqref{WSh}  is closely linked to the theory of the Hankel transform	(see, for instance, \cite{Watson}, p. 398) and its  representation depends on whether $x_1$ is less than, equal to, or greater than $x_2$.
In the case $m{\,=\,}3$, the infinite integral \eqref{WSh} was independently computed by Bailey (\cite{Bailey}, p. 45, Eq. (7.1)) and Rice (\cite{Rice}, p. 60, Eq. (2.6)), each employing distinct approaches. 

All the integrals presented in this article are evaluated using the methodology developed in the handbook by Marichev \cite{Marichev1983}. This work is notable for its dual focus: it provides both a rigorous theoretical framework and a practical, algorithmic approach to evaluating complex integrals. By systematically applying this methodology, one can derive a vast array of integral formulas involving higher transcendental functions, including the Weber-Schafheitlin integrals considered herein.

 \section{Pearson random walk }
 
 Pearson's original papers \cite{PEARSON1,PEARSON2} gave birth to the entire field of random walk theory, and the Pearson random walk remains the simplest and most elegant model for isotropic diffusion in two dimensions. Its mathematical beauty lies in the deep connection between probability theory, Bessel functions, and the heat equation.

Pearson states that migration is a statistical problem, yet it also demands careful consideration of geographical and ecological factors, including food and shelter availability, which are specific to each species. The role of mathematics is not to provide a definitive model for every real-world case, but rather to construct an idealized framework that offers valuable insights when applied judiciously.

Pearson's idealized model rested on the following assumptions:
\begin{enumerate}
	\item Breeding sites and food resources are uniformly distributed across the region, with no preferential movement along rivers or forest paths.
\item The species scattering from a point of origin occurs symmetrically in all directions. The characteristic distance an individual moves between habitats is termed a "flight"\,, and an organism may undergo $n$ such flights—either from its origin to a breeding site or between successive breeding grounds if multiple reproductive events occur. This "flight" is distinct from a "flitter"\,, which denotes short, localized back-and-forth movements for feeding or mating.
\end{enumerate}

  {\bf  Problems:} {\it Given some numbers of individuals departing from a single point (idealized as a geometric point), what is their spatial distribution after $n$ random flights? This constitutes the core question, which Pearson termed the Fundamental Problem of Random Migration.
    The next stage involves distributing these points across a region with an arbitrary boundary and tracking the resulting spatial distribution on both sides of that boundary over multiple breeding seasons.  These are the Secondary Problems of Migration.}

  The answer involves integrating over all possible angles and this angular integration is precisely what produces the product of Bessel functions.

This problem has a long interesting story.
 From the inception of the Educational Times in 1847 and its later Mathematical Questions series, numerous problems have pertained to probability theory, with particular emphasis on random walks. The earliest documented mathematical treatment appears to be that of M. W. Crofton \cite{crofton}, who examined a traveler's motion along a river, formalizing what is now understood as a linear random walk.
 Major Ross in \cite{Ross} considered the problem as being one of the infiltration of mosquitoes into cleared areas  when he studied the laws of the spread of malaria. The canonical formulation, however, is due to Karl Pearson \cite{PEARSON1,PEARSON2}, who in 1905 inquired in Nature about the radial distance distribution after $n$ unit steps with random orientations. Lord Rayleigh's immediate response indicated prior resolution of this problem in the context of acoustic wave propagation \cite{Rayleigh0}. Shortly thereafter, Jan Cornelis Kluyver   advanced a comprehensive solution in terms of Bessel-function integrals \cite{Kluyver1905}. Pearson and his collaborator in \cite{Pearson1906} subsequently applied Kluyver's formalism to a migration study, employing graphical techniques and manual power-series approximations due to the absence of computational resources and the limited state of special-function theory at the time. Notably, all of these contributions were restricted to planar walks; Watson later (see \cite{Watson}) extended the framework to arbitrary dimensions. 
 The first model with a hyperbolic governing equation was proposed by Sydney Goldstein \cite{Goldstein1951}, who studied a one-dimensional random walk with velocity reversals driven by a Poisson process. This was later elaborated by Mark Katz \cite{Katz1}, Enzo Orsingher, and others \cite{Orsingher1985,ors1,ors2,ors3,garra,Kelbert,garra21}, with further generalizations given in subsequent articles. For additional historical detail, the reviews by Jacques Dutka \cite{dutka},  Subrahmanyan Chandrasekhar \cite{chandra}  are recommended.
Nowadays, ``The Ultimate Univariate Probability Distribution Explorer''~\cite{blog} 
  contains about 30000 formulas for 500 most known probability distributions and their 60 properties  
  (including the Pearson distribution and abstract distribution). 

 \section{Probabilistic interpretation of Weber-Schafheitlin integral}
 
 Now we consider mathematical problem.
 Suppose that the $N$-dimensional random vectors $\{\mathbf{r}^j\}$, $\mathbf{r}^j=(r^j_1,...,r^j_N)$ with tips which are uniformly distributed  on a $N$-dimensional sphere centered in origin. Their sum
$$
 	\mathbf{R} = \mathbf{r}^1 + \mathbf{r}^2 + \dots + \mathbf{r}^n=\sum\limits_{j=1}^n \mathbf{r}^j.
$$
 has a radially symmetric distribution, meaning its probability density function  depends only on 
  $R = |\mathbf{R}|$ only.  For independent, identically distributed  $\{\mathbf{r}^j\}$, this probability density function  is conveniently obtained via the characteristic function. 
Characteristic function for the random vector $\mathbf{r}^j$ is:  
$$
 \varphi_{\mathbf{r}^j}(\boldsymbol{\xi})= \mathbb{E}[e^{i\langle \boldsymbol{\xi}, \mathbf{r}^j \rangle}] = \int\limits_{\mathbb{R}^N} e^{i\langle \boldsymbol{\xi}, \boldsymbol{\rho}\rangle} f_{\mathbf{r}_j}(\boldsymbol{\rho}) \, d\boldsymbol{\rho},
 $$
 where
$$ 
\boldsymbol{\xi}=(\xi_1,...,\xi_N),\qquad \boldsymbol{\rho}=(\rho_1,...,\rho_N),\qquad \langle \boldsymbol{\xi}, \boldsymbol{\rho} \rangle=\sum\limits_{i=1}^N \xi_i\rho_i.
$$  

Since, in our case the random vector $\{\mathbf{r}^j\}$ is  uniformly distributed in direction on an $N${-}dimensional sphere, then its
probability density function $f_{\mathbf{r}^j}(\boldsymbol{\rho})$  is  supported on an $N$-dimensional  sphere   and 
$$
\varphi_{\mathbf{r}^j}(\boldsymbol{\xi})=\frac{1}{\omega_N(|\mathbf{r}^j|)}\int\limits_{S^N(|\mathbf{r}^j|)} e^{i\langle \boldsymbol{\xi}, \boldsymbol{\rho}\rangle}  \, dS_N,
$$
 where $S^N(|\mathbf{r}^j|)$ is a $N$-dimensional  sphere with radius $|\mathbf{r}^j|$,  
 $dS_{N}$ is a surface measure, $\omega_N(|\mathbf{r}^j|){\,=\,}\frac{2\pi^{\frac{N}{2}}|\mathbf{r}^j|^{N-1}}{\Gamma\left(\frac{N}{2}\right)}$ is the total surface area of the unit sphere in $\mathbb{R}^N$.  Under the integral  there is a plane wave function $e^{i\langle \boldsymbol{\xi}, \boldsymbol{\rho}\rangle} $.
 
Moving on to the integration over the unit sphere, we obtain
 $$
 \varphi_{\mathbf{r}^j}(\boldsymbol{\xi})=\frac{1}{\omega_N(1)}\int\limits_{S^N(1)} e^{i|\mathbf{r}^j|\langle \boldsymbol{\xi}, \boldsymbol{\rho}\rangle}  \, dS_N.
 $$
 
  Due to the 
 rotational symmetry (isotropy) of the uniform distribution  on the sphere, the  
characteristic function depend only on the Euclidean norm $|\boldsymbol{\xi}|{\,=\,}\sqrt{\xi_1^2+...+\xi_N^2}$ of the vector $\boldsymbol{\xi}$. We can rotate the coordinate system such that the vector $\boldsymbol{\xi}$ points entirely along the first axis: $\boldsymbol{\xi}\to (|\boldsymbol{\xi}|,0,...,0)$. Then, the scalar product $\langle \boldsymbol{\xi}, \boldsymbol{\rho}\rangle$ simplified to:
$$
\langle \boldsymbol{\xi}, \boldsymbol{\rho}\rangle\to |\boldsymbol{\xi}| \rho_1.
$$
 
 Let $\theta \in [0,\pi]$ be defined by:
 $\rho_1{\,=\,}\cos \theta.$
  Geometrically, $\theta$ is the polar angle (colatitude) from the positive $\rho_1$-axis.
  The remaining $(N-1)$ coordinates $(\rho_2, \ldots, \rho_N)$ of unit sphere must satisfy:
 $$
 \rho_2^2 + \cdots + \rho_N^2 = 1 - \rho_1^2 = 1 - \cos^2 \theta = \sin^2 \theta.
 $$
  So the vector $(\rho_2, \ldots, \rho_N)$ lies on a sphere of radius $\sin \theta$ in $\mathbb{R}^{N-1}$ and  the full surface element is:
  $dS_N = (\sin \theta)^{N-2} \, d\theta \, dS_{N-1}$.
  
   Therefore, we get 
  $$
  \varphi_{\mathbf{r}^j}(\boldsymbol{\xi})=\frac{1}{\omega_N(1)}\int\limits_{S^N(1)} e^{i|\mathbf{r}^j|\cdot|\boldsymbol{\xi}| \rho_1}  \, dS_N=\frac{1}{\omega_N(1)}\int\limits_0^\pi e^{i|\mathbf{r}^j|\cdot|\boldsymbol{\xi}|\cos\theta} (\sin \theta)^{N-2} d\theta \int\limits_{S^{N-1}(1)} dS_{N-1}=
  $$
  $$
  =\frac{\omega_{N-1}(1)}{\omega_N(1)}\int\limits_0^\pi e^{i|\mathbf{r}^j|\cdot|\boldsymbol{\xi}|\cos\theta} (\sin \theta)^{N-2} d\theta=\frac{\Gamma\left(\frac{N}{2}\right)}{\sqrt{\pi}\Gamma\left(\frac{N-1}{2}\right)}\int\limits_0^\pi e^{i|\mathbf{r}^j|\cdot|\boldsymbol{\xi}|\cos\theta} (\sin \theta)^{N-2} d\theta.
  $$
Recalling the standard integral representation for the Bessel function of the first kind \cite{Watson}:
$$J_\nu(z) = \frac{\left(\frac{z}{2}\right)^\nu}{\sqrt{\pi} \, \Gamma\left(\nu + \frac{1}{2}\right)} \int\limits_0^\pi e^{iz \cos \theta} (\sin \theta)^{2\nu} \, d\theta,$$
valid for $\operatorname{Re}\,(\nu) > -\frac{1}{2}$, we obtain
$$
 \varphi_{\mathbf{r}^j}(\boldsymbol{\xi})=
 \frac{\Gamma\left(\frac{N}{2}\right)}{\sqrt{\pi}\Gamma\left(\frac{N-1}{2}\right)}\cdot\frac{\sqrt{\pi}\Gamma\left(\frac{N-1}{2}\right)}{\left(\frac{|\mathbf{r}^j|\cdot|\xi|}{2}\right)^{\frac{N}{2}-1}} J_{\frac{N}{2}-1}(|\mathbf{r}^j|\cdot|\boldsymbol{\xi}|)=\frac{2^{\frac{N}{2}-1}\Gamma\left(\frac{N}{2}\right)}{(|\mathbf{r}^j|\cdot|\xi|)^{\frac{N}{2}-1}}J_{\frac{N}{2}-1}(|\mathbf{r}^j|\cdot|\boldsymbol{\xi}|).
$$ 

%Let us notice that the general  fundamental identity for integral by the sphere $S^{N}(1)$ is
%$$
%	\int\limits_{S^{N}(1)} g(\langle y, x\rangle) \, dS_N 
%	= \omega_{N}(1) \int\limits_{-1}^{1} (1 - t^2)^{\frac{N-3}{2}} g(|y| t) \, dt,
%$$
%then 

For the sum $\mathbf{R} = \sum\limits_{j=1}^n \mathbf{r}_j$, the characteristic function is:
$$
\varphi_{\mathbf{R}}(\boldsymbol{\xi}) = \mathbb{E}[e^{i\langle \boldsymbol{\xi}, \mathbf{R} \rangle}] = \prod\limits_{j=1}^n \varphi_{\mathbf{r}^j}(\boldsymbol{\xi})= |\boldsymbol{\xi}|^{n\left( 1-\frac{N}{2}\right)}\left(2^{\frac{N}{2}-1}\Gamma\left(\frac{N}{2}\right)\right)^n\prod\limits_{j=1}^n |\mathbf{r}^j|^{ 1-\frac{N}{2}}\cdot J_{\frac{N}{2}-1}(|\mathbf{r}^j|\cdot|\boldsymbol{\xi}|).
$$

 	The probability density function $p(\mathbf{R})$ is the inverse Fourier transform of the characteristic function:
$$
p(\mathbf{R}) = \frac{ 2^{N\left(\frac{n}{2}-1\right) -n}\Gamma^n\left(\frac{N}{2}\right)}{\pi^N} \int\limits_{\mathbb{R}^N} |\boldsymbol{\xi}|^{n\left( 1-\frac{N}{2}\right)}e^{-i\langle \boldsymbol{\xi}, \mathbf{R} \rangle} \prod_{j=1}^n |\mathbf{r}^j|^{ 1-\frac{N}{2}}\cdot J_{\frac{N}{2}-1}(|\mathbf{r}^j|\cdot|\boldsymbol{\xi}|)   d\boldsymbol{\xi}.
$$

Converting points from Cartesian   into spherical coordinates $\boldsymbol{\xi}=\rho \Theta$ we get
$$
p(\mathbf{R}) = \frac{ 2^{N\left(\frac{n}{2}-1\right) -n}\Gamma^n\left(\frac{N}{2}\right)}{\pi^N} \int\limits_0^\infty   \rho^{n\left( 1-\frac{N}{2}\right)} \prod_{j=1}^n |\mathbf{r}^j|^{ 1-\frac{N}{2}}\cdot J_{\frac{N}{2}-1}(|\mathbf{r}^j|\cdot \rho) \rho^{N-1}d\rho\int\limits_{S^N(1)} e^{-i\rho\langle \Theta, \mathbf{R} \rangle} dS_N.
$$

%$$
%\int\limits_{S^N(1)} e^{-i\rho \langle \Theta, \mathbf{R} \rangle} \, dS_N 
%= \frac{2\pi^{\frac{N}{2}}}{\Gamma\left(\frac{N}{2}\right)} 
%\frac{J_{\frac{N-1}{2}}(\rho)}{(\rho/2)^{\frac{N-1}{2}}},
%$$

Calculating integral by sphere as earlier by the formula 
\begin{equation*}
	\int\limits_{S^N(1)} e^{-i\rho\langle \Theta, \mathbf{R} \rangle} dS_N=\frac{(2\pi)^{\frac{N}{2}}}{(\rho |\mathbf{R}|)^{\frac{N}{2}-1}}J_{\frac{N}{2}-1}(\rho |\mathbf{R}|),
\end{equation*}
 we get probability density function in the form
\begin{multline}\label{IntDen}
 p(\mathbf{R})=\\
  = \frac{ 2^{n\left(\frac{N}{2}-1\right)-\frac{N}{2}}\Gamma^n\left(\frac{N}{2}\right)}{\pi^\frac{N}{2}|\mathbf{R}|^{ \frac{N}{2}-1}} \int\limits_0^\infty   \rho^{n\left( 1-\frac{N}{2}\right)+\frac{N}{2}}  J_{\frac{N}{2}-1}(|\mathbf{R}|\cdot \rho) \prod_{j=1}^n |\mathbf{r}^j|^{ 1-\frac{N}{2}}\cdot J_{\frac{N}{2}-1}(|\mathbf{r}^j|\cdot \rho)  d\rho.
\end{multline}
 
 Taking into account spherical symmetry the cumulative distribution function has the form
 $$
 F(\mathbf{R})=\omega_N(1)\int\limits_0^{|\mathbf{R}|} p(\mathbf{R})|_{|\mathbf{R}|=r} r^{N-1}dr,
 $$
 then
\begin{multline*}
 F(\mathbf{R})=\omega_N(1)\frac{ 2^{n\left(\frac{N}{2}-1\right)-\frac{N}{2}}\Gamma^n\left(\frac{N}{2}\right)}{\pi^\frac{N}{2}}\times\\
 \times \int\limits_0^\infty   \rho^{n\left( 1-\frac{N}{2}\right)+\frac{N}{2}}  \left( \int\limits_0^{|\mathbf{R}|}r^{ \frac{N}{2}}J_{\frac{N}{2}-1}(\rho r)dr\right)  \prod_{j=1}^n |\mathbf{r}^j|^{ 1-\frac{N}{2}}\cdot J_{\frac{N}{2}-1}(|\mathbf{r}^j|\cdot \rho)  d\rho
\end{multline*}
 or, calculating integral
 \begin{multline}\label{IntPr}
 	 F(\mathbf{R})=2^{n\left(\frac{N}{2}-1\right)-\frac{N}{2}+1}\Gamma^{n-1}\left(\frac{N}{2}\right)|\mathbf{R}|^{\frac{N}{2}}\times\\
\times \int\limits_0^\infty   \rho^{n\left( 1-\frac{N}{2}\right)+\frac{N}{2}-1}   J_{\frac{N}{2}}(|\mathbf{R}| \rho )  \prod_{j=1}^n |\mathbf{r}^j|^{ 1-\frac{N}{2}}\cdot J_{\frac{N}{2}-1}(|\mathbf{r}^j|\cdot \rho)  d\rho.
\end{multline}

	In \cite{Kuznetsov}  a simulation and analytical study of a multidimensional random walk of multiple agents, where each agent may turn by any angle at each step was  presented. The analytical model gives the probability of finding an agent within a radius $r$ at a given time using generalization of integral \eqref{IntPr}, with parameters for motion intensity and environment.

  \section{Calculation of the Weber--Schafheitlin integral}
 
{\bf Generalized Weber--Schafheitlin integral and special cases.}
 Integrals \eqref{IntDen} and \eqref{IntPr} has the form of the Weber--Schafheitlin integral \eqref{WSh}.
 
 In \cite{Srivastava} the next result was obtained.

 \begin{teo} 
 Let $x_1, \ldots, x_m$ be positive real numbers such that
 $$
 x_m > x_1 + \cdots + x_{m-1}, \quad m \geq 2.
 $$
 Also let $M{\,=\,}\mu+\nu_1+...+\nu_m$. Then
 \begin{multline}\label{eq:2.8}
 	\int\limits_{0}^{\infty} t^{\mu -1} \prod_{j = 1}^{m} J_{\nu_j}(x_j t)\, dt=  \\
 	= \frac{2^{\mu - 1} x_1^{\nu_1} \cdots x_{m-1}^{\nu_{m-1}} x_m^{\nu_m - M}
 		\Gamma\left(\frac{M}{2}\right)}
 	{\Gamma(\nu_1 + 1)\cdots\Gamma(\nu_{m-1} + 1)
 		\Gamma\left(\nu_m - \frac{M}{2} + 1\right)} \times\\
 	\times F_c^{(m-1)}
 	\left[\frac{M}{2}, \frac{M}{2} - \nu_m;
 	\nu_1 + 1, \ldots, \nu_{m-1} + 1;
 	\frac{x_1^2}{x_m^2}, \ldots, \frac{x_{m-1}^2}{x_m^2}\right],
 \end{multline}
 provided that
$$
 	\operatorname{Re}(1 + \nu_1 + \cdots + \nu_m) > \operatorname{Re}(1 - \mu) > -\frac{m}{2}.
 $$
 \end{teo}

In \eqref{eq:2.8} $F_c^{(m-1)}$ is Lauricella's function. 
{\bf Lauricella's hypergeometric function} $F_c^{(n)}$ of $(n)$ variables $x{\,=\,}x_1,...,x_n$
with complex parameters $\alpha$, $\beta$,  $\gamma{\,=\,}(\gamma_1, \ldots, \gamma_n)$ is defined by
(cf. \cite{Lauricella}, p. 113)
\begin{multline}\label{laur}
	F_c^{(n)}(\alpha, \beta,\gamma;x)=F_c^{(n)}(\alpha, \beta; \gamma_1, \ldots, \gamma_n; x_1, \ldots, x_n)=\\
	= \sum_{j_1, \ldots, j_n = 0}^{\infty}
	\frac{(\alpha)_{j_1 + \cdots + j_n}(\beta)_{j_1 + \cdots + j_n}}
	{(\gamma_1)_{j_1}\cdots(\gamma_n)_{m_n}}
	\frac{x_1^{j_1}\cdots x_n^{j_n}}{j_1!\cdots j_n!}
\end{multline}
with, as usual,
$$
(\lambda)_j = \frac{\Gamma(\lambda + j)}{\Gamma(\lambda)} =
\begin{cases}
	1, & \text{if } j = 0; \\
	\lambda(\lambda + 1)\cdots(\lambda + j - 1), & \forall j \in \{1, 2, 3, \ldots\}.
\end{cases}
$$
The series \eqref{laur} converges in the domain
$$
D=\{(x_1,...,x_n)\in\mathbb{C}^n:|x_1|^{1/2} + \cdots + |x_n|^{1/2} < 1\}.
$$
However, outside this domain, the function $F_c^{(n)}$ can be continued to the region
$
\sum\limits_{i=1}^{n-1} \sqrt{|x_i|}{\,+\,}1{\,<\,}\sqrt{|x_n|}
$
by means of the formula
\begin{multline*}
F_c^{(n)}(\alpha, \beta,\gamma;x)=	F_C^{(n)}(\alpha,\beta;\gamma_1,\dots,\gamma_n;x_1,\dots,x_n)= \\
	= \frac{\Gamma(\gamma_n)\Gamma(\beta-\alpha)}{\Gamma(\gamma_n-\alpha)\Gamma(\beta)} (-x_n)^{-\alpha}
	F_C^{(n)}\!\left(\alpha, 1+\alpha-\gamma_n; \gamma_1,\dots,\gamma_{n-1}, 1+\alpha-\beta; 
	\frac{x_1}{x_n}, \dots, \frac{x_{n-1}}{x_n}, \frac{1}{x_n}\right)+ \\
	+ \frac{\Gamma(\gamma_n)\Gamma(\alpha-\beta)}{\Gamma(\gamma_n-\beta)\Gamma(\alpha)} (-x_n)^{-\beta}
	F_C^{(n)}\!\left(\beta, 1+\beta-\gamma_n; \gamma_1,\dots,\gamma_{n-1}, 1+\beta-\alpha; 
	\frac{x_1}{x_n}, \dots, \frac{x_{n-1}}{x_n}, \frac{1}{x_n}\right).
\end{multline*}

Let $f(x){\,=\,}f(x_1, \dots, x_n)$ be a function of $n$ complex variables. Define the differential operators
$$
\theta_k = x_k \frac{\partial}{\partial x_k}, \qquad 
\theta = \sum_{k=1}^n \theta_k.
$$
The Lauricella function $F_C^{(n)}(\alpha, \beta, \gamma; x)$ is a solution to the following system of $n$ second-order partial differential equations (see \cite{Goto}):
$$
\left[ \theta_k (\theta_k + \gamma_k - 1) - x_k (\theta + \alpha)(\theta + \beta) \right] f(x) = 0, \qquad k = 1, \dots, n.
$$
The system generated by them is called Lauricella’s hypergeometric system $E_C (\alpha, \beta, \gamma)$ of differential equations.

As was shown in \cite{Srivastava}
formula \eqref{eq:2.8} simplifies when
$$
	\mu = \nu_m - \nu_1 - \cdots - \nu_{m-1},
$$
the Lauricella series on the right-hand side reduces to 1, and  the following elegant result take place (cf. \cite{Srivastava}, p. 47, Eq. (10.4))
$$
\int\limits_{0}^{\infty} t^{\nu_m - \nu_1 - \cdots - \nu_{m-1} - 1}
\prod_{j = 1}^{m} J_{\nu_j}(x_j t)\, dt  = \frac{2^{\nu_m - \nu_1 - \cdots - \nu_{m-1} - 1} x_1^{\nu_1} \cdots x_{m-1}^{\nu_{m-1}} \Gamma(\nu_m)}
{x_m^{\nu_m} \Gamma(\nu_1 + 1)\cdots\Gamma(\nu_{m-1} + 1)},
\label{eq:3.2}
$$
provided that $x_m > x_1 + \cdots + x_{m-1}$. 
  
{\bf Generalized Weber--Schafheitlin integral.} 
	There are four classical kinds of Bessel functions: 
\begin{enumerate}
	\item $J_\nu(x)$ is the Bessel function of the first kind 
	 is defined by the series
	$$
		J_\nu(x) = \sum\limits_{m=0}^\infty \frac{(-1)^m}{m!\, \Gamma(m+\nu+1)} {\left({\frac{x}{2}}\right)}^{2m+\nu}.
	$$
	For non-integer $\alpha$ the functions $J_\nu (x)$ and $J_{-\nu} (x)$ are linearly independent.
	If $\nu$ is integer  the following relationship is valid:
	$
	J_{-\nu}(x){\,=\,}(-1)^{\nu} J_{\nu}(x).
	$
	This function is finite at $x = 0$ for $\nu \ge 0$. 
	\item The Bessel function of the second kind (Weber/Neumann function), denoted by $Y_\nu(x)$, is related to $J_\nu(x)$ for non-integer $\nu$ by the formula:
	$$
	Y_\nu(x) = \frac{J_\nu(x) \cos(\nu\pi) - J_{-\nu}(x)}{\sin(\nu\pi)}.
$$
	In the case of integer order $n$, the function $Y_\nu(x)$ is defined by taking the limit as the non-integer $\nu$ tends to $n$:
	$ Y_{n}(x){\,=\,}\lim\limits_{\nu \to n}Y_{\nu }(x).
	$
This function is singular at $x = 0$.
	\item The Hankel Function of the first kind (gives outgoing wave):	
	$$
	H_\nu^{(1)}(x) = J_\nu(x) + iY_\nu(x).
	$$
	\item The Hankel function of the second kind (gives incoming wave):	
	$$
	H_\nu^{(2)}(x) = J_\nu(x) - iY_\nu(x).
	$$
\end{enumerate}

The functions $J_\nu(x)$, $Y_\nu(x)$, $H_\nu^{(1)}(x)$, and $H_\nu^{(2)}(x)$ are all solutions to the Bessel equation of order $\nu$:
$$
	x^2 \frac{d^2 y}{dx^2} + x \frac{dy}{dx} + (x^2 - \nu^2) y = 0.
$$

	There are two more common Bessel functions:
		\begin{enumerate}
		\item 	The modified Bessel function of the first kind, denoted by $I_\nu(x)$, is defined by the series:
				$$
			I_\nu(x)=\frac{x^{\alpha}}{(i x)^{\alpha}} J_{\alpha}(i x) = \sum_{m=0}^{\infty} \frac{1}{m! \, \Gamma(m+\nu+1)} \left(\frac{x}{2}\right)^{2m+\nu}.
		$$
		 This function is real-valued for real $x$ and is finite at $x = 0$ for $\nu \ge 0$. 
			\item The modified Bessel function of the second kind, denoted by $K_\nu(x)$, is defined in terms of $I_\nu(x)$ and $I_{-\nu}(x)$:
			$$
			K_\nu(x) = \frac{\pi}{2} \frac{I_{-\nu}(x) - I_{\nu}(x)}{\sin(\nu\pi)}.
		$$
		For integer order $n$, this is defined by the limiting form:
		$K_n(x) = \lim\limits_{\nu \to n} K_\nu(x)$. The function $K_\nu(x)$ is singular at $x = 0$ and decays exponentially as $x \to \infty$.
	\end{enumerate}
	
	The functions $I_\nu(x)$  and $K_\nu(x)$ are all solutions to the modified Bessel equation of order~$\nu$:
	$$
		x^2 \frac{d^2 y}{dx^2} + x \frac{dy}{dx} - (x^2 + \nu^2) y = 0.
	$$

	So   get six total in the ``standard family'' of Bessel functions.

			Spherical Bessel functions of arbitrary order $\nu$
			 arise when solving the Helmholtz equation in spherical coordinates. They are closely related to the ordinary Bessel functions of order $\nu + \frac{1}{2}$:
			 			 \begin{enumerate}
			 	\item Spherical Bessel function of the first kind of order $\nu$, denoted by $j_\nu(x)$:
			 	$$
			 		j_\nu(x) = \sqrt{\frac{\pi}{2x}} \, J_{\nu+\frac{1}{2}}(x).
			 	$$
			 	This function is finite at $x = 0$ for $\nu\ge -\frac{1}{2}$.
			 		\item  Spherical Bessel function of the second kind of order $\nu$ (also called the spherical Neumann function), denoted by $y_\nu(x)$:
			 	$$
			 		y_\nu(x) = \sqrt{\frac{\pi}{2x}} \, Y_{\nu+\frac{1}{2}}(x).
			 	$$
			 	This function is singular at $x = 0$.
			 		\item  Spherical Hankel function of the first kind of order $\nu$, denoted by $h_\nu^{(1)}(x)$:
			 	$$
			 		h_\nu^{(1)}(x) = j_\nu(x) + i y_\nu(x) = \sqrt{\frac{\pi}{2x}} \, H_{\nu+\frac{1}{2}}^{(1)}(x).
			 	$$
			 	This represents an outgoing spherical wave.
			 		\item Spherical Hankel function of the second kind of order $\nu$, denoted by $h_\nu^{(2)}(x)$:
			 	$$
			 		h_\nu^{(2)}(x) = j_\nu(x) - i y_\nu(x) = \sqrt{\frac{\pi}{2x}} \, H_{\nu+\frac{1}{2}}^{(2)}(x).
			 	$$
			 	This represents an incoming spherical wave.
			 \end{enumerate}

	The functions $j_\nu(x)$, $y_\nu(x)$, $h_\nu^{(1)}(x)$, and $h_\nu^{(2)}(x)$ are all solutions to the spherical Bessel equation of order $\nu$:	
			$$
				x^2 \frac{d^2 y}{dx^2} + 2x \frac{dy}{dx} + \left[x^2 - \nu(\nu+1)\right] y = 0.
		$$
	
	This gives a total of ten core functions in the extended Bessel family.

We have two more function related to Bessel function of the first kind, namely, the confluent hypergeometric limit function ${}_0F_1$ which is related to the Bessel function of the first kind $J_\nu$ by the following identity:
		$$
					{}_0F_1\left(; \nu; z\right) = - \,\frac{z \, \Gamma(\nu)}{(i \sqrt{z})^{\nu+1}} \, J_{\nu-1}\left(2i \sqrt{z}\right),
		$$
				where  $\nu \notin \{0, -1, -2, \dots\}$ and the regularized confluent hypergeometric limit function ${}_0\widetilde{F}_1$ is related to the Bessel function of the first kind $J_\nu$ by:
$$
										{}_0\widetilde{F}_1\left(; \nu; z\right) = -\,\frac{z}{(i \sqrt{z})^{\nu+1}}  J_{\nu-1}\left(2i \sqrt{z}\right).
$$

Let us consider Weber--Schafheitlin integral \eqref{WSh}.
 We are able to construct much more "super general"\, integrals, including \eqref{WSh} as particular cases, which can be expanded of finite sums of Weber--Schafheitlin integrals. Indeed, for non integer parameter $\nu$ we have relations
	\begin{align*}\label{rel}
		Y_\nu(z) &= J_\nu(z) \cot(\pi \nu) - \csc(\pi \nu) \, J_{-\nu}(z), \\
		H_\nu^{(1)}(z) &= \left(1 + i \cot(\pi \nu)\right) J_\nu(z) - i \csc(\pi \nu) \, J_{-\nu}(z), \\
		H_\nu^{(2)}(z) &= \left(1 - i \cot(\pi \nu)\right) J_\nu(z) + i \csc(\pi \nu) \, J_{-\nu}(z), \\
		I_\nu(z) &= (i z)^{-\nu} z^\nu \, J_\nu(i z), \\
		K_\nu(z) &= \frac{\pi}{2} \csc(\pi \nu) \left[ (i z)^\nu z^{-\nu} J_{-\nu}(i z) - (i z)^{-\nu} z^\nu J_\nu(i z) \right], \\
		j_\nu(z) &= \sqrt{\frac{\pi}{2z}} \, J_{\nu+\frac{1}{2}}(z), \\
		y_\nu(z) &= -\sqrt{\frac{\pi}{2z}} \sec(\pi \nu) \, J_{-\frac{1}{2}-\nu}(z) - \sqrt{\frac{\pi}{2z}} \tan(\pi \nu) \, J_{\frac{1}{2}+\nu}(z), \\
		h_\nu^{(1)}(z) &= -\sqrt{\frac{\pi}{2z}} \, i \sec(\pi \nu) \, J_{-\frac{1}{2}-\nu}(z) + \sqrt{\frac{\pi}{2z}} \left(1 - i \tan(\pi \nu)\right) J_{\frac{1}{2}+\nu}(z), \\
		h_\nu^{(2)}(z) &= i \sec(\pi \nu) \sqrt{\frac{\pi}{2z}} \, J_{-\nu-\frac{1}{2}}(z) + (1 + i \tan(\pi \nu)) \sqrt{\frac{\pi}{2z}} \, J_{\nu+\frac{1}{2}}(z), \\
		{}_0F_1\left(; \nu; z\right) &= -(i \sqrt{z})^{-1-\nu} \, z \, \Gamma(\nu) \, J_{-1+\nu}\left(2i \sqrt{z}\right),\\
		{}_0\widetilde{F}_1\left(; \nu; z\right) &= -(i \sqrt{z})^{-1-\nu} \, z \, J_{\nu-1}\left(2i \sqrt{z}\right).
	\end{align*}
	These formulas demonstrate that eleven functions on the left-hand side can be expressed as linear combinations of products with function $J_\mu(\xi)$.
 Consequently, we can now define the maximal general Weber–Schafheitlin integral as:
$$
		\begin{aligned}
			&\text{GWS}\left[\beta, \left\{n_j\right\}_{j=1}^{12}\right]= \\
			&= \int\limits_0^\infty t^{\beta-1}
			\left( \prod_{j=1}^{n_1} J_{\nu_{j,1}}(x_{j,1} t) \right)
			\left( \prod_{j=1}^{n_2} Y_{\nu_{j,2}}(x_{j,2} t) \right)
			\left( \prod_{j=1}^{n_3} H_{\nu_{j,3}}^{(1)}(x_{j,3} t) \right)
			\left( \prod_{j=1}^{n_4} H_{\nu_{j,4}}^{(2)}(x_{j,4} t) \right)\times \\
			&\quad \times
			\left( \prod_{j=1}^{n_5} I_{\nu_{j,5}}(x_{j,5} t) \right)
			\left( \prod_{j=1}^{n_6} K_{\nu_{j,6}}(x_{j,6} t) \right)
			\left( \prod_{j=1}^{n_7} j_{\nu_{j,7}}(x_{j,7} t) \right)
			\left( \prod_{j=1}^{n_8} y_{\nu_{j,8}}(x_{j,8} t) \right)
			\left( \prod_{j=1}^{n_9} h_{\nu_{j,9}}^{(1)}(x_{j,9} t) \right)\times \\
			&\quad \times
						\left( \prod_{j=1}^{n_{10}} h_{\nu_{j,10}}^{(2)}(x_{j,10} t) \right)
			\left( \prod_{j=1}^{n_{11}} {}_0F_1\left(; \nu_{j,11}; x_{j,11} t^2 \right) \right)\left(\prod_{j=1}^{n_{12}} {}_0\widetilde{F}_1\left(; \nu_{j,12}; x_{j,12} t^2 \right)\right) 
				\, dt,
		\end{aligned}
$$
 where the parameters satisfy appropriate convergence conditions.
	
This integral  for non integer indexes $\nu$  is equal to a finite sum of Weber--Schafheitlin integrals \eqref{WSh} with corresponding parameters and coefficients. For example,
for the parameter configuration
$\{n_1, 0, 0, 0, n_5, 0, n_7, 0, 0, 0, 0, n_{12}\}$,
	the Generalized Weber--Schafheitlin integral reduces to:
	$$
	\operatorname{WS}[\mu;\overline{\nu};\overline{x}] =\int\limits_{0}^{\infty} t^{\mu-1}\prod\limits_{j=1}^m J_{\nu_j}(x_jt)\,dt,
$$
		\begin{equation}\label{WSexample}
		\begin{aligned}
			&\text{WS}\left[\mu;\overline{\nu};\overline{x}\right]= \\
			&= C \int\limits_0^\infty t^{\mu - 1}
			\prod_{j=1}^{n_1} J_{\nu_{j,1}}(x_{j,1} t)
			\prod_{j=1}^{n_5} J_{\nu_{j,5}}(i x_{j,5} t)
			\prod_{j=1}^{n_7} J_{\nu_{j,7} + 1/2}(x_{j,7} t)
			\prod_{j=1}^{n_{12}} J_{\nu_{j,12} - 1}\left(2i \sqrt{x_{j,12}} \, t \right)
			\, dt,
		\end{aligned}
	\end{equation}
		where
	\begin{align*}
		m &= n_1 + n_5 + n_7 + n_{12},\\
		\overline{\nu}&=(\nu_{1,1},...,\nu_{n_1,1},\nu_{1,5},...,\nu_{n_5,5},\nu_{1,7}+ 1/2,...,\nu_{n_7,7}+ 1/2,
		\nu_{1,12}- 1,...,\nu_{n_{12},12}- 1),\\
		\overline{x}&=(x_{1,1}, \ldots, x_{n_1,1},\; i x_{1,5}, \ldots, i x_{n_5,5},\; x_{1,7}, \ldots, x_{n_7,7},\; 2i\sqrt{x_{1,12}}, \ldots, 2i\sqrt{x_{n_{12},12}} ),\\
		C &= \exp\left( \frac{\pi i}{2} \left( \sum_{j=1}^{n_{12}} (1 - \nu_{j,12}) - \sum_{j=1}^{n_5} \nu_{j,5} \right) \right)
		\left( \frac{\pi}{2} \right)^{n_7/2}
		\left( \prod_{j=1}^{n_7} x_{j,7} \right)^{-1/2}
		\left( \prod_{j=1}^{n_{12}} x_{j,12}^{(1 - \nu_{j,12})/2} \right), \\
		\mu &= \beta - \frac{n_7}{2} + \sum_{j=1}^{n_{12}} (1 - \nu_{j,12}).
	\end{align*}
which is the classical Weber--Schafheitlin integral \eqref{WSh}.
	A necessary condition for convergence at infinity of \eqref{WSexample} is that all entries of the vector $\overline{x}$ are real.

  \section{Particular cases of  $m$}

	For the $m{\,=\,}1$
	$$
		\int\limits_{0}^{\infty} t^{\mu-1} J_{\nu_1}(x_1t) \, dt 
		= \frac{2^{\mu-1} \, \Gamma\!\left( \frac{\nu_1+\mu}{2}\right)}{x_1^\mu\Gamma\!\left( \frac{\nu_1-\mu}{2} +1\right)}, 
		\qquad \operatorname{Re}\,(\mu+ \nu_1) > 0, \;\; \operatorname{Re}\,(\mu) < \frac{3}{2},\qquad x_1>0.
	$$

In particular, when $\mu{\,=\,}1$  (see 10.22.41 from \cite{DLMF}):
$$
\int\limits_{0}^{\infty} J_{\nu_1}(t)dt=1,\qquad \operatorname{Re}(\nu_1)>-1.
$$

In the  case when $\mu=1$, $\nu_1=0$ we obtain	
$$
	\int\limits_{0}^{\infty} J_0(xt)dt=\left\{\begin{array}{cc}
		-\dfrac{1}{x}, & x<0; \\ [0.3 cm]
		\dfrac{1}{x}, &  x>0. 
	\end{array} \right.
$$

Setting $m{\,=\,}2$ in \eqref{WSh}  and consulting \cite{WS}, \cite{DLMF} \S10.22.56–57, we obtain a formula that evaluates the   integral \eqref{WSh} as a
Gauss hypergeometric function.

General formula here is
\begin{multline}\label{TwoJ}
	\int\limits_{0}^{\infty} t^{\mu-1} J_{\nu_1}(x_1 t) J_{\nu_2}(x_2 t) \, dt=\\
	=
	\begin{cases}
		\dfrac{2^{\mu-1} x_2^{-\mu-\nu_1} x_1^{\nu_1} \, \Gamma\!\left(\frac{\mu+\nu_1+\nu_2}{2}\right)}
		{\Gamma\!\left(\frac{\nu_2-\mu-\nu_1}{2}+1\right) \Gamma(\nu_1+1)}
		\, {}_2F_1\!\left(\dfrac{\mu+\nu_1+\nu_2}{2}, \dfrac{\mu+\nu_1-\nu_2}{2}; \nu_1+1; \dfrac{x_1^2}{x_2^2}\right),
		& x_2 > x_1; \\[1.2em]
		\dfrac{2^{\mu-1} x_1^{-\mu-\nu_2} x_2^{\nu_2} \, \Gamma\!\left(\frac{\mu+\nu_1+\nu_2}{2}\right)}
		{\Gamma\!\left(\frac{\nu_1-\mu-\nu_2}{2}+1\right) \Gamma(\nu_2+1)}
		\, {}_2F_1\!\left(\dfrac{\mu+\nu_1+\nu_2}{2}, \dfrac{\mu-\nu_1+\nu_2}{2}; \nu_2+1; \dfrac{x_2^2}{x_1^2}\right),
		& x_1 > x_2; \\[1.2em]
		\dfrac{2^{\mu-1} x_1^{-\mu} \, \Gamma\!\left(\frac{\mu+\nu_1+\nu_2}{2}\right) \Gamma(1-\mu)}
		{\Gamma\!\left(\frac{\nu_1-\mu-\nu_2}{2}+1\right)
			\Gamma\!\left(\frac{\nu_2-\mu-\nu_1}{2}+1\right)
			\Gamma\!\left(\frac{\nu_1+\nu_2-\mu}{2}+1\right)},
		& x_1 = x_2, \\[1.2em]
		\text{does not converge in classical case}, & \text{otherwise}
	\end{cases}
\end{multline}
with  $x_1 > 0$, $x_2 > 0$, $\operatorname{Re}\,(\mu + \nu_1 + \nu_2) > 0$, $\operatorname{Re}(\mu) < 2$ for $x_1\neq x_2$ and $\operatorname{Re}(\mu) < 1$ for $x_1=x_2$.

About 30 particular cases of \eqref{TwoJ} for positive $x_1$ and $x_2$ one can find at chapters 2.12.31-2.12.32 of handbook \cite{IR2}, below we present some of such cases, which also available at \cite{DLMF}  (see 10.22.58, 10.22.55,  10.22.62, 10.22.63, 10.22.64, from \cite{DLMF}).

When $\nu_1{\,=\,}\nu_2$, we get formula
\begin{multline*}
\int\limits_0^\infty t^{\mu-1} J_{\nu_1}(x_1 t) J_{\nu_1}(x_2 t)\,dt
=\\
= 
\frac{
	2^{\mu-1} (x_1x_2)^{\nu_1} (x_1^2+x_2^2)^{-\frac{\mu}{2}-\nu_1} \Gamma\!\left(\frac{\mu}{2}+\nu_1\right)
}{
	\Gamma\!\left(1-\frac{\mu}{2}\right)
}
\,
{}_2F_1\!\left(
\frac{\mu+2\nu_1}{4},
\frac{2+\mu+2\nu_1}{4},
1+\nu_1,
\frac{4x_1^2x_2^2}{(x_1^2+x_2^2)^2}
\right)
\end{multline*}
provided that
$x_1 \ne x_2$, $\operatorname{Re}\,(\mu+2\nu_1)>0$,  $\operatorname{Re}\,(\mu)<2$.

In particular, if $\nu_2{\,=\,}\nu_1$ and $x_2{\,=\,}x_1$, we obtain
$$
\int\limits_0^\infty t^{\mu-1} J_{\nu_1}^2(x_1 t) \, dt
= \frac{2^{\mu-1}  \Gamma(1-\mu) \Gamma\!\left(\frac{\mu}{2}+\nu_1\right)}
{\Gamma^2\!\left(1-\frac{\mu}{2}\right) \Gamma\!\left(1-\frac{\mu}{2}+\nu_1\right)}\,x_1^{-\mu},
\quad \operatorname{Re}\,(\mu+2\nu_1)>0,\; \operatorname{Re}\,(\mu)<1.
$$

When $\mu=\nu_1-\nu_2+2$  we obtain
$$
\int\limits_0^\infty t^{\nu_1-\nu_2+1} J_{\nu_1}(x_1 t) J_{\nu_2}(x_2 t) \, dt
= 
\begin{cases}
	0, & 0 < x_2 < x_1; \\[6pt]
	\dfrac{2^{\nu_1-\nu_2+1} x_1^{\nu_1} (x_2^2-x_1^2)^{\nu_2-\nu_1-1}}{ \Gamma(\nu_2-\nu_1)x_2^{\nu_2}}, & x_1<x_2,
\end{cases}
$$
under conditions $-1 < \operatorname{Re}\,(\nu_1) < \operatorname{Re}\,(\nu_2)$.

	For the special case $\mu{\,=\,}1$, $\nu_1{\,=\,}\nu_2{-}1$,  formula \eqref{TwoJ} can be used to describe the probability that a particle exits a circle of a given radius. This is based on the following special case (see \cite{Weinstein1948} and \cite{Watson}, p. 461):
$$
\int\limits_0^\infty J_{\nu_2-1}(x_1t)J_{\nu_2}(x_2t)dt=\left\{ \begin{array}{ll}
	\dfrac{x_1^{\nu_2-1}}{x_2^{\nu_2}}, & \mbox{$0<x_1<x_2$};\\ [0.3 cm]
	\dfrac{1}{2x_2}, & \mbox{$x_1=x_2$};\\ [0.3 cm]
	0 & \mbox{for $x_1>x_2$}.\end{array} \right.
$$

We have another interesting 
formula when $\mu{\,=\,}1$ and $\nu_1{\,=\,}\nu+2n+1$, $\nu_2{\,=\,}\nu$
$$
\int\limits_0^\infty J_{\nu+2n+1}(x_1 t) J_\nu(x_2 t) \, dt
=
\begin{cases}
	\dfrac{x_2^\nu \Gamma(\nu+n+1)}{x_1^{\nu+1} \Gamma(\nu+1)n!}
	\, {}_2{F}_1\!\left(-n, \nu+n+1; \nu+1; \dfrac{x_2^2}{x_1^2}\right), & 0 < x_2 < x_1; \\[8pt]
	\dfrac{(-1)^n}{2x_1}, & x_1 = x_2 > 0; \\[4pt]
	0, & x_1<x_2,
\end{cases}
$$
under conditions $\operatorname{Re}\,(\nu)>-n-1,\; n\in\mathbb{Z},\; n\ge 0.$

One of the simplest but impressive case  when  $\mu{\,=\,}1$ is $\nu_1{\,=\,}\nu_2=0$:
$$
\int\limits_0^\infty J_0(x_1 t) J_0(x_2 t) \, dt
= \frac{2 }{\pi |x_1|}K\!\left(\frac{x_2^2}{x_1^2}\right),
\qquad \frac{x_2^2}{x_1^2} < 1.
$$
Here $K(z){\,=\,}\int\limits_{0}^{\pi/2} \frac{dt}{\sqrt{1 - z \sin^2 t}}$ 
gives the complete elliptic integral of the first kind when
 $\left|{\rm Arg}(1 - z) \right|{\,<\,}\pi$.

Now turn to the case $\mu{\,=\,}2$.
Although the integral \eqref{TwoJ} diverges in the classical sense for $\nu_1{\,=\,}\nu_2$ and $\operatorname{Re}(\mu) < 1$,
it can nevertheless be evaluated when understood in the sense of distributions.
For example, when $\mu=2$, $\nu_1=\nu_2$, we have
$$
\int\limits_{0}^{\infty} t J_{\nu_1}(x_1 t) J_{\nu_1}(x_2 t) \, dt = \frac{1}{\sqrt{x_1 x_2}} \delta(x_1-x_2)
$$
is used to prove that the inversion of the Hankel transform
$
F_\nu(\xi){\,=\,} \int\limits_{0}^{\infty} f(x) J_\nu(x \xi) \, x\, dx
$
for $\nu{\,>\,}{-}1/2$,  is given by the identical self-reciprocal integral formula
$
f(x){\,=\,}\int\limits_{0}^{\infty} F_\nu(\xi) J_\nu(x \xi) \, \xi \, d\xi.
$ 

Interesting cases appear when $\mu{\,=\,}0$ and $x_1{\,=\,}x_2$:  	
$$
\int\limits_0^\infty \frac{1}{t} J_{\nu_1}(x_1 t) J_{\nu_2}(x_1 t) \, dt
= \frac{\Gamma\!\left(\frac{\nu_1+\nu_2}{2}\right)}
{2\,\Gamma\!\left(\frac{\nu_2-\nu_1}{2}+1\right)
	\Gamma\!\left(\frac{\nu_1-\nu_2}{2}+1\right)
	\Gamma\!\left(\frac{\nu_1+\nu_2}{2}+1\right)},
\quad \operatorname{Re}\,(\nu_1+\nu_2)>0.
$$
In particular, the last formula gives orthogonality property of Bessel functions. Namely, for $\mu{\,=\,}0$, $\nu_1{\,=\,}\nu+2\ell+1$, $\nu_2{\,=\,}\nu+2m+1$, $x_1{\,=\,}x_2{\,=\,}1$, $m,\ell{\,\in\,}\mathbb{N}_0$ (see \S 10.22.55 from \cite{DLMF})):
$$
\int\limits_0^\infty \frac{1}{t} J_{\nu+2\ell+1}(t) J_{\nu+2m+1}(t) \, dt
= \frac{\delta_{\ell m}}{2(2\ell+\nu+1)},
$$
where $\delta_{\ell m}$ is a Kronecker delta.
When $\mu{\,=\,}0$ and $\nu_1{\,=\,}\nu_2$, $x_1{\,=\,}x_2$
$$
\int\limits_0^\infty \frac{1}{t} J_{\nu_1}(x_1 t)^2 \, dt
= \frac{1}{2\nu_1}, \quad \operatorname{Re}\,(\nu_1)>0.
$$

		For the $m{\,=\,}3$, in the general case  we get
	the following formula (see \cite{Bailey}, p.~45, Eq.~(7.1)); \cite{Rice}, p.~60, Eq.~(2.6) and \cite{WS1}), where below integral was also represented through the Kamp\'{e} de F\'{e}riet function
	\begin{multline}\label{ThreeJ}
		\int\limits_{0}^{\infty} t^{\mu-1} J_{\nu_1}(x_1 t) J_{\nu_2}(x_2 t) J_{\nu_3}(x_3 t) \, dt=\\
		=\mathbf{C}\cdot F_4\left(\frac{\mu +\nu_1+\nu_2 -\nu_3}{2},\frac{\mu+\nu_1 +\nu_2 +\nu_3}{2} ,\nu_1 +1,\nu_2
		+1,\frac{x_1^2}{x_3^2},\frac{x_2^2}{x_3^2}\right)=\\
				=\mathbf{C}\cdot
F^{2,0,0}_{0,1,1}
\left(\begin{matrix}
	\frac{\mu+\nu_1+\nu_2-\nu_3}{2};& \frac{\mu+\nu_1+\nu_2+\nu_3}{2};&- ;&-;\,\\[0.2em]
	-;&	\nu_1+1;& \nu_2+1;&
\end{matrix}
\frac{x_1^2}{x_3^2}, \frac{x_2^2}{x_3^2}
\right)=\\
=\mathbf{C}\cdot F_C^{(2)}\left( \frac{\mu+\nu_1+\nu_2-\nu_3}{2},\frac{\mu+\nu_1+\nu_2+\nu_3}{2},\nu_1+1,\nu_2+1,\frac{x_1^2}{x_3^2}, \frac{x_2^2}{x_3^2}
\right),
\end{multline}
where
$$
\mathbf{C}=\frac{2^{\mu-1} x_1^{\nu_1}x_2^{\nu_2} x_3^{-\mu-\nu_1-\nu_2} \,
	\Gamma\!\left(\frac{\mu+\nu_1+\nu_2+\nu_3}{2}\right)}
{\Gamma(\nu_1+1) \Gamma(\nu_2+1) \,
	\Gamma\!\left(1-\frac{\mu+\nu_1+\nu_2-\nu_3}{2}\right)}
$$
			with $x_1,x_2,x_3 \in \mathbb{R}$,   $x_1,x_2,x_3 > 0$, $x_3>x_1+x_2$,
	$\operatorname{Re}(\mu+\nu_1+\nu_2+\nu_3) > 0$, $\operatorname{Re}(\mu) < \frac{5}{2}$.

In \eqref{ThreeJ} $F_4$ is the {\bf Appell function}, $F_{2,0,0}^{0,1,1}$ is the {\bf Kamp\'{e} de F\'{e}riet function} and $F_C^{(2)}$ is {\bf Lauricella function}:
	\begin{multline*}
			F_4(a,b,c_1,c_2,z,w)
				=F^{2,0,0}_{0,1,1}\left({\begin{matrix}a,b;-;-;\\
					-;c_1;c_2;
			\end{matrix}}z,w\right)=F_C^{(2)}(a,b,c_1,c_2,z,w)=\\
					=\sum _{m=0}^{\infty } \sum _{n=0}^{\infty } \frac{ (a)_{m+n} (b)_{m+n}}{\left(c_1\right)_m
		\left(c_2\right)_n}\frac{z^m w^n}{m! n!},\qquad \sqrt{|z| }+\sqrt{|w| }<1.
	\end{multline*}
The Appell function $F_4$ and the Lauricella function $F_C^{(n)}$  are particular cases of the more general Kamp\'e de F\'eriet function $F_{P,Q,S}^{A,B,C}$, which  is a special function that generalizes a hypergeometric function to two variables:
$$
F^{A,B,C}_{P,Q,S}
\left({\begin{matrix}a_{1},\cdots,a_{A};b_{1},\cdots,b_{B}; c_{1},\cdots ,c_{C};\\
		p_{1},\cdots p_{P};q_{1},\cdots, q_Q; s_1,\cdots, s_S;
\end{matrix}}z,w\right)=\sum _{m=0}^{\infty }\sum _{n=0}^{\infty }\frac{\prod\limits_{j=1}^A(a_{j})_{m+n}
	\prod\limits_{j=1}^B(b_{j})_{m}\prod\limits_{j=1}^C(c_{j})_{n} }{\prod\limits_{j=1}^P(p_{j})_{m+n}
	\prod\limits_{j=1}^Q(q_{j})_{m}\prod\limits_{j=1}^S(s_{j})_{n} } \frac {z^mw^n}{m!n!},
$$
where
\begin{multline*}
	(A+B\leq P+Q\land A+C\leq P+S\land | z| <\infty \land | w| <\infty )\lor\\
	\lor(A+B=P+Q+1\land A+C<P+S+1\land | z| <1\land | w| <\infty
)\lor\\
\lor(A+B<P+Q+1\land A+C=P+S+1\land | z| <\infty \land | w| <1)\lor\\
\lor (A+B=P+Q+1\land A+C=P+S+1\land | z| <1\land | w| <1\land
A\leq P)\lor\\
\lor \left(A+B=P+Q+1\land A+C=P+S+1\land | w| ^{\frac{1}{A-P}}+| z| ^{\frac{1}{A-P}}<1\land A>P\right).
\end{multline*}
	
The Appell function $F_4$ can be analytically continued outside the domain
$
\sqrt{|z|} + \sqrt{|w|} < 1
$
to the domain
$
\sqrt{|z|} + 1 < \sqrt{|w|}
$
by applying the following formula:
$$
\begin{aligned}
	F_4(&a,b;c_1,c_2;z,w)= \\
	&= \frac{\Gamma(c_2)\Gamma(b-a)}{\Gamma(c_2-a)\Gamma(b)} (-w)^{-a}
	F_4\!\left(a, 1+a-c_2; c_1, 1+a-b; \frac{z}{w}, \frac{1}{w}\right)+ \\
	&\quad + \frac{\Gamma(c_2)\Gamma(a-b)}{\Gamma(c_2-b)\Gamma(a)} (-w)^{-b}
	F_4\!\left(b, 1+b-c_2; c_1, 1+b-a; \frac{z}{w}, \frac{1}{w}\right).
\end{aligned}
$$

Impressive case of formula 	\eqref{ThreeJ} with $\nu_1{\,=\,}\nu_2{\,=\,}\nu_3{\,=\,}\nu$  and $\mu{\,=\,}=2-\nu$ was given in \cite{IR2} p. 204, formulas 2.12.42.14-15 in the form 
		\begin{equation}\label{ThreeJ1}
\int\limits_0^{\infty}t^{1-\nu}J_\nu(at)J_\nu(bt)J_\nu(ct)dt=
\frac{2^{1-\nu}}{\sqrt{\pi}(abc)^\nu\Gamma\left(\nu+\frac{1}{2}\right) }\varDelta^{2\nu-1}, 
\end{equation}
where $ a,b,c{\,>\,}0$, $(b{+}c{-}a)(a{+}c{-}b)(a{+}b{-}c){\,>\,} 0$,  
$\varDelta{\,=\,}\frac{1}{4}\sqrt{(b{+}c{-}a)(a{+}c{-}b)(a{+}b{-}c)(a{+}b{+}c)}$,
$|a{-}b|{\,<\,}c{\,<\,}a{\,+\,}b$, 
$\operatorname{Re}\,\nu{\,>\,}{-}\frac{1}{2}$ and $\varDelta$ is the area of a triangle whose sides are equal to $a$, $b$ and $c$.

The most important thing that is formula \eqref{ThreeJ1} allows us to express the solution of the Cauchy problem of the form
$$
\left\{
\begin{array}{l}
\left(\dfrac{\partial^2}{\partial x^2}+\dfrac{\gamma}{x}\dfrac{\partial}{\partial x}\right)   u(x,y)=\left(\dfrac{\partial^2}{\partial y^2}+\dfrac{\gamma}{y}\dfrac{\partial}{\partial y}\right) u(x,y), \\ [0.3 cm]
u(x,0)=f(x),\\
\dfrac{\partial}{\partial y}u(x,y)\biggr|_{y=0}=0.
\end{array}
\right.
$$
 Namely,
  \begin{equation*}\label{FitSv}
 u(x,y)=\left( \frac{2}{xy}\right) ^{\frac{\gamma-1}{2}} \Gamma\left(\frac{\gamma+1}{2}\right) 
 	\int\limits_{|x-y|}^{x+y}z^{\frac{5-3\gamma}{2}} f(z) dz\int\limits_0^{\infty} J_{\frac{\gamma-1}{2}}(xt) J_{\frac{\gamma-1}{2}}(yt) J_{\frac{\gamma-1}{2}}(zt)
 	t^{\frac{3-\gamma}{2}} \, dt.
 \end{equation*}
 This solution $u(x,y){\,=\,}\,^\gamma T_x^yf(x)$ is the generalized translation 
 (see \cite{sitnikshishkinaElsevier}, p. 146, Eq. (3.155)).

 About 45 particular cases of \eqref{ThreeJ1} for positive $x_1$, $x_2$, $x_3$ 
 one can find at chapter 2.12.42 of handbook [33], below we present just some of such cases, which also available at \cite{IR2}   (see formulas 10.22.71, 10.22.72 from \cite{IR2}). When $\mu{\,=\,}2{-}\nu_1$, $\nu_2{\,=\,}\nu_3$ we get
 $$
\int\limits_0^{\infty}t^{1-\nu_1}J_{\nu_1}(x_1t)J_{\nu_2}(x_2t)J_{\nu_2}(x_3t)dt =\frac{ (x_2x_3)^{\nu_1 -1} \sin ^{\nu_1 -\frac{1}{2}}(\varphi )}{\sqrt{2 \pi
	}x_1^{\nu_1 }} P_{\nu_2 -\frac{1}{2}}^{\frac{1}{2}-\nu_1}\cos (\varphi),
$$
where
$$
\operatorname{Re}(\nu_1 )>-\frac{1}{2}\land \operatorname{Re}(\nu_2 )>-1\land | x_2-x_3| <x_1<x_2+x_3\land \cos (\varphi )=\frac{x_2^2+x_3^2-x_1^2}{2 x_2 x_3},
 $$
 and
 $$
\int\limits_0^{\infty}t^{1-\nu_1}J_{\nu_1}(x_1t)J_{\nu_2}(x_2t)J_{\nu_2}(x_3t)dt =
\frac{\sqrt{2} e^{i \pi  \left(\nu_1 -\frac{1}{2}\right)}   \sin (\pi  (\nu_1 -\nu_2 )) 
	\sinh ^{\nu_1
		-\frac{1}{2}}(\chi ) }{\pi ^{3/2}x_1^{\nu_1 }(x_2x_3 )^{1-\nu_1}}Q_{\nu_2 -\frac{1}{2}}^{\frac{1}{2}-\nu_1 }(\cosh (\chi )),
$$
where
$$
\operatorname{Re}(\nu_1 )>-\frac{1}{2}\land
\operatorname{Re}(\nu )>-1\land x_1>x_2+x_3\land \cosh (\chi )=\frac{x_1^2-x_2^2-x_3^2}{2 x_2x_3}.
$$
Here $P_n^m(x)$ is a Legendre function of the first kind and
$Q_n^m(x)$ is a Legendre function of the second kind.

Additional infinite integrals over the product of three Bessel functions (including modified Bessel functions) are given in Gervois and Navelet  \cite{GervoisNavelet1984,GervoisNavelet1985a,GervoisNavelet1985b,GervoisNavelet1986a,GervoisNavelet1986b}.

	Using symmetry and permutations $(\nu_1, x_1) \leftrightarrow (\nu_2, x_2) \leftrightarrow (\nu_3, x_3)$, 
	we can write analogical formulas for $x_1 > x_2 + x_3$ and $x_2 > x_1 + x_3$.  
	For other regions like $|x_2-x_3| < x_1 < x_2+x_3$ 
	the value of integral still was not described in all situations. 

\section{Conclusion}

In this paper, we have presented a comprehensive study of the Weber–Schafheitlin integral and its profound connections to random walk theory. We have shown that the probability distributions for multidimensional random walks naturally lead to improper integrals involving products of Bessel functions, precisely of the Weber–Schafheitlin type. This connection provides a powerful analytical framework for evaluating radial probability density functions and cumulative distribution functions in isotropic random flight models.

\end{document}